\documentclass{article}
\usepackage{amsmath, amssymb, amsthm, epsf, graphicx, amscd, amsfonts, amscd, multirow}

\usepackage[pagewise, displaymath, mathlines]{lineno}

\newtheorem{theorem}{Theorem}[section]
\newtheorem{conjecture}[theorem]{Conjecture}
\newtheorem{lemma}[theorem]{Lemma}

\newtheorem{corollary}[theorem]{Corollary}

\theoremstyle{definition}

\def\bb #1{ {\mathbb #1} }
\def\c #1{ {\mathcal #1} }

\def\b #1{ {\bf #1} }

\DeclareMathOperator{\lcm}{lcm}
\DeclareMathOperator{\ord}{ord}

\begin{document}
\title{Partial zeta functions, partial exponential sums, \\ and $p$-adic estimates }

\author{Noah Bertram, Xiantao Deng, C. Douglas Haessig, Yan Li}

\date{\today}
\maketitle

\abstract{Partial zeta functions of algebraic varieties over finite fields generalize the classical zeta function by allowing each variable to be defined over a possibly different extension field of a fixed finite field. Due to this extra variation their rationality is surprising, and even simple examples are delicate to compute. For instance, we give a detailed description of the partial zeta function of an affine curve where the number of unit poles varies, a property different from classical zeta functions. On the other hand, they do retain some properties similar to the classical case. To this end, we give Chevalley-Warning type bounds for partial zeta functions and $L$-functions associated to partial exponential sums.}


\section{Introduction}

Let $\bb F_q$ be the finite field with $q = p^a$ elements, $p$ a prime number. Fix $\b d := (d_1, \ldots, d_n) \in \bb Z_{\geq 1}^n$, and set $d := \lcm(d_1, \ldots, d_n)$. Throughout the paper, we will write $q_i := q^{d_i}$. Let $X$ be an affine variety defined over $\bb F_{q^d}$, defined as the zero locus of $F_1, \ldots, F_r \in \bb F_{q^d}[x_1, \ldots, x_n]$. Associated to the sequence
\[
m \geq 1: \qquad N_m(\b d) := \# \{ x := (x_1, \ldots, x_n) \in \bb F_{q_1^m} \times \cdots \times \bb F_{q_n^m} \mid F_1(x) = \cdots = F_r(x) = 0 \}
\]
is the \emph{partial zeta function} of $X$,
\[
Z(X / \bb F_{q}, \b d, T) := \exp\left( \sum_{m \geq 1} N_m(\b d) \frac{T^m}{m} \right).
\]
Wan introduced partial zeta functions in \cite{MR1803946} and proved their rationality in \cite{MR2027782}:
\begin{equation}\label{E: zeta}
Z(X / \bb F_{q}, \b d, T) = \frac{\prod_{i=1}^R (1 - \alpha_i T)}{\prod_{j = 1}^S (1 - \beta_j T)} \in \bb Q(T),
\end{equation}
with $\alpha_i$ and $\beta_j$ Weil $q$-integers. Since each variable lies in a possibly different finite extension field of $\bb F_q$, even simple examples can be complicated, and very different from their classical counterpart. Further, their rationality is surprising, even in the most simple cases. This point is illustrated in Sections \ref{S: hyperplane} and \ref{S: affine curve} where we take a detailed look at a couple of examples. In particular, the partial zeta function of the affine curve $\c C_n$ defined by $y = x^n$ takes the form
\begin{equation}\label{E: Cn example}
Z(\c C_n / \bb F_{q}, \b d, T) = \frac{1}{1-T} \cdot \prod_{k \mid \varphi(n)}\left( \frac{\Phi_k(T)}{\Phi_k(q^cT)} \right)^{a_k},
\end{equation}
where $\Phi_k$ is the $k$-th cyclotomic polynomial, $\varphi$ is Euler's totient function, $c = \gcd(d_1, d_2)$, and
\[
a_k := \frac{1}{\varphi(n)} \sum_{i = 1}^{\varphi(n)}  \gcd(n, M_i) \zeta_k^{i} \in \bb Z,
\]
where $\zeta_k$ is a primitive $k$-th root of unity and $M_i :=  (q_1^i - 1)/(q^{ci}-1)$. An interesting feature of this example is that, for fixed $n$ and $q$, the number of unit poles varies between zero and one as $d_1$ and $d_2$ vary. This does not occur for classical zeta functions, where affine curves never have unit poles.

We note that while Fu and Wan \cite{MR2030375} have given a cohomological description of these zeta functions, rationality is still not immediate since their trace formula involves roots of unity, which translates into the zeta function having factors with roots of unity as exponents. Rationality is obtained by a further Galois theoretic argument. We wondered whether this Galois argument is an artifact of their method or intrinsic to these types of zeta functions, and indeed, it appears to be intrinsic. In the example of the affine curve $\c C_n$, we examine its partial zeta function using an ad hoc method, and as you can see above, the exponents $a_k$ are exponential sums.

Little is known about the $p$-adic behavior of the zeros and poles of these types of zeta functions. A first step toward such an understanding is the following theorem, proven in Section \ref{S: partial zeta}. To state the result we need to define some notation. For a nonnegative integer $b$, write $b = a_0 + a_1 p + \cdots + a_r p^r$ with $0 \leq a_i \leq p-1$. Set $\sigma_p(b) := a_0 + \cdots + a_r$. Extend this to $u = (u_1, \ldots, u_n) \in \bb Z_{\geq 0}^n$ by defining $|\sigma_p(u)| := \sigma_p(u_1) + \cdots + \sigma_p(u_n)$. Using multi-index notation $x^u := x_1^{u_1} \cdots x_n^{u_n}$, for a polynomial $F(x) = \sum a_u x^u \in \bb F_q[x_1, \ldots, x_n]$ define the $p$-weight $w_p(F) := \max_u |\sigma_p(u)|$. 

\begin{theorem}\label{T: main zeta}
Let $X$ be an affine variety defined as the zero locus of $F_1, \ldots, F_r \in \bb F_{q^d}[x_1, \ldots, x_n]$. Then $p^\omega \mid N_1(\b d)$, where
\[
\omega := \left\lceil a \cdot \frac{(d_1 + \cdots + d_n) - d \sum_{i=1}^r w_p(F_i)}{\max_i w_p(F_i)}  \right\rceil.
\]
\end{theorem}

A similar theorem was proven in \cite[Theorem 1.4]{MR1803946}.

\begin{corollary}\label{C: zeta ord}
$\ord_q \alpha_i$ and $\ord_q \beta_j$ from (\ref{E: zeta}) are bounded below by
\[
\frac{(d_1 + \cdots + d_n) - d \sum_{i=1}^r w_p(F_i)}{\max_i w_p(F_i)}.
\] 
\end{corollary}

We also prove an analogue of Theorem \ref{T: main zeta} for $L$-functions of partial exponential sums, which are defined as follows. Let $f(x) \in \bb F_{q^d}[x_1, \ldots, x_n]$, $\psi$ a nontrivial additive character on $\bb F_{q^d}$, and $\chi_1, \ldots, \chi_r$ multiplicative characters on $\bb F_{q^{d_1}}^\times, \ldots, \bb F_{q^{d_n}}^\times$, respectively. Define the mixed partial character sum
\begin{align*}
S({\bf \chi}, {\bf d}, f) := \sum \chi_1(x_1) \cdots \chi_r(x_r) \psi( f(x) ),
\end{align*}
where the sum runs over $x_i \in \bb F_{q^{d_i}}^\times$ for $i = 1, \ldots, r$ and $x_i \in \bb F_{q^{d_i}}$  for  $i = r+1, \ldots, n$. To give a $p$-adic estimate for $S({\bf \chi}, {\bf d}, f)$, we embed it into the $p$-adic numbers. 

Let $\bb Q_p$ be the $p$-adic numbers, and denote by $\bb Q_{q^d}$ the unramified extension of $\bb Q_p$ of degree $ad$ (recall, $q = p^a$). We normalize the valuation such that $\ord_q(q) = 1$. Denote by $T_i$ the set of solutions of $z^{q_i} = z$ in $\bb Q_{q^{d_i}} \subseteq \bb Q_{q^d}$, and set $T_i^\times := T_i \setminus \{0\}$. Define
\[
\c T := T_1^\times \times \cdots \times T_r^\times \times T_{r+1} \times \cdots \times T_n.
\]
There exist $e_i \in \{1, \ldots, q_i - 1\}$ for $i = 1, \ldots, r$ such that  
\[
S(\chi, \b d, f) = \sum_{x \in \c T} x_1^{e_1} \cdots x_r^{e_r} \psi(f(\bar x)) \in \bb Q_{q^d}(\zeta_p),
\]
where $\bar x_i$ denotes the image of $x_i$ in the residue field $\bb F_{q^{d_i}} \subseteq \bb F_{q^d}$.

\begin{theorem}
\[
\ord_q S(\chi, \b d, f) \geq \frac{1}{w_p(f)} \left( \sum_{i = 1}^n d_i -  \frac{1}{a(p-1)} \sum_{i=1}^r \sigma_p(e_i)  \right).
\]
When there is no twist (i.e. $r = 0$) the estimate becomes
\[
\ord_q S(\b d, f) \geq (d_1 + \cdots + d_n) / w_p(f).
\]
\end{theorem}

{\bf Acknowledgments.} Xiantao Deng thanks Yang Zhang for helpful discussions. Doug Haessig thanks Prof. Daqing Wan for valuable discussions, and introducing him to Deng and Li.

\section{Example: Affine hyperplane}\label{S: hyperplane}

Let $\b d := (d_1, \ldots, d_n) \in \bb Z_{\geq 1}^n$. Observe that affine space $\bb A^n$ has the partial zeta function
\[
Z(\bb A^n / \bb F_q, \b d, T) = \frac{1}{1 - q^{d_1 + \cdots + d_n} T}.
\]

Next, consider the affine hyperplane $H$ defined by $a_1 x_1 + \cdots. + a_n x_n = 0$, with $a_i \in \bb F_q^*$. 

\begin{theorem}
\begin{equation}\label{E: hyperplane}
Z(H / \bb F_q, \b d, T) = \frac{1}{1 - q^e T}, \qquad \text{where} \qquad e := \sum_{k = 2}^n  (-1)^{k} \sum_{1 \leq i_1 < \cdots < i_k \leq n} \gcd(d_{i_1}, \ldots, d_{i_k}).
\end{equation}
\end{theorem}

The proof of the formula for $e$ requires the following lemma. 

\begin{lemma}\label{L: Fq results}
\phantom{}
\begin{enumerate}
\item  $\bb F_{q^{d_1}} \cap (\bb F_{q^{d_2}} + \cdots + \bb F_{q^{d_n}}) = \bb F_{q^{d_1}} \cap \bb F_{q^{d_2}} + \cdots + \bb F_{q^{d_1}} \cap \bb F_{q^{d_n}}$.
\item The dimension of $\bb F_{q^{d_1}} + \cdots + \bb F_{q^{d_n}}$ as an $\bb F_q$-vector space is $\sum_{k=1}^n (-1)^{k+1} \sum_{1 \leq i_1 < \cdots < i_k \leq n} \gcd( d_{i_1}, \ldots, d_{i_k})$.
\end{enumerate}
\end{lemma}

\begin{proof}
For the first statement, we need only show
\[
\bb F_{q^{d_1}} \cap (\bb F_{q^{d_2}} + \cdots + \bb F_{q^{d_n}}) \subseteq \bb F_{q^{d_1}} \cap \bb F_{q^{d_2}} + \cdots + \bb F_{q^{d_1}} \cap \bb F_{q^{d_n}}
\]
since it is clear the right hand side is contained in the left. Set $m := \lcm(d_1, \ldots, d_n)$. Assume for the moment that $\ord_p d_1 \geq \ord_p d_i$ for every $i$, and thus $p \nmid m / d_1$. Write $\bb F_{q^m}^* = \langle \zeta \rangle = \{ \zeta^i \mid 1 \leq i \leq q^m - 1\}$, where $\zeta$ is a primitive $q^m-1$ root of unity. Observe that $\bb F_{q^{d_i}}^* = \langle \zeta^{a_i} \rangle$, where $a_i := (q^m-1) / (q^{d_i} - 1)$. Let $\eta \in \bb F_{q^{d_1}} \cap (\bb F_{q^{d_2}} + \cdots + \bb F_{q^{d_n}})$. Since $\eta \in \bb F_{q^{d_2}} + \cdots + \bb F_{q^{d_n}}$ we may write $\eta = \sum_{i=2}^n b_i \zeta^{a_i t_i}$ for some $b_i \in \bb F_p$ and $t_i \in \bb Z$. Next, since $\eta \in \bb F_{q^{d_1}}$, then for all $j \in \bb Z$ we have
\[
\eta = \eta^{q^{j d_1}} = \sum_{i=2}^n b_i \zeta^{a_i t_i q^{j d_1}}.
\]
Thus, setting $u := m / d_1$ then $c \eta = \sum_{i=2}^n b_i \sum_{j = 0}^{u-1} \zeta^{a_i t_i q^{j d_1}}$. Notice that
\[
\left( \sum_{j=0}^{u-1} \zeta^{a_i t_i q^{j d_1}} \right)^{q^{d_1}} = \sum_{j=0}^{u-1} \zeta^{a_i t_i q^{j d_1}},
\]
and so $ \sum_{j=0}^{u-1} \zeta^{a_i t_i q^{j d_1}} \in \bb F_{q^{d_1}} \cap \bb F_{q^{d_i}}$. Since $p \nmid u$, it follows that
\[
\eta = \sum_{i=1}^n \frac{b_i}{u} \sum_{j=0}^{c-1} \zeta^{a_i t_i q^{j d_1}} \in \bb F_{q^{d_1}} \cap \bb F_{q^{d_2}} + \cdots. + \bb F_{q^{d_1}} \cap \bb F_{q^{d_n}}
\]
as desired. For this argument, we made the assumption that $\ord_p d_1 \geq \ord_p d_i$ for every $i$. We may remove this restriction after proving the second statement of the lemma by induction as follows.

It is well-known that $\bb F_{q^{a}} \cap \bb F_{q^b} = \bb F_{q^{\gcd(a, b)}}$, and so when $n = 2$ we have
\begin{align*}
\dim_{\bb F_q}( \bb F_{q^{d_1}} + \bb F_{q^{d_2}}) &= \dim_{\bb F_q} \bb F_{q^{d_1}}  + \dim_{\bb F_q} \bb F_{q^{d_2}}  - \dim_{\bb F_q}( \bb F_{q^{d_1}} \cap \bb F_{q^{d_2}}) \\
&= \dim_{\bb F_q} \bb F_{q^{d_1}}  + \dim_{\bb F_q} \bb F_{q^{d_2}}  - \dim_{\bb F_q}( \bb F_{q^{\gcd(d_1, d_2)}}) \\
&= d_1 + d_2 - \gcd(d_1, d_2).
\end{align*}
Suppose now that $n > 2$. Reorder the $d_i$ so that $\ord_p d_1 \geq \ord_p d_i$. Set $r_i := \gcd(d_1, d_i)$. Then
\begin{align*}
\dim_{\bb F_q}(\bb F_{q^{d_1}} + \cdots + \bb F_{q^{d_n}}) &= \dim_{\bb F_q} \bb F_{q^{d_1}} + \dim_{\bb F_q}(\bb F_{q^{d_2}} + \cdots + \bb F_{q^{d_n}}) - \dim_{\bb F_q}\left( \bb F_{q^{d_1}} \cap (\bb F_{q^{d_2}} + \cdots + \bb F_{q^{d_n}}) \right) \\
&= \dim_{\bb F_q} \bb F_{q^{d_1}} + \dim_{\bb F_q}(\bb F_{q^{d_2}} + \cdots + \bb F_{q^{d_n}}) - \dim_{\bb F_q}(\bb F_{q^{r_2}} + \cdots + \bb F_{q^{r_n}})
\end{align*}
Hence,
\begin{align*}
\dim_{\bb F_q}(\bb F_{q^{d_1}} + \cdots + \bb F_{q^{d_n}}) &= d_1 + \sum_{k = 1}^{n-1} (-1)^{k+1} \sum_{2 \leq i_1 < \cdots < i_k \leq n} \gcd(d_{i_1}, \ldots, d_{i_k}) -  \sum_{k = 1}^{n-1} (-1)^{k+1} \sum_{2 \leq i_1 < \cdots < i_k \leq n} \gcd(r_{i_1}, \ldots, r_{i_k}) \\ 
&= d_1 + \sum_{k = 1}^{n-1} (-1)^{k+1} \sum_{2 \leq i_1 < \cdots < i_k \leq n} \gcd(d_{i_1}, \ldots, d_{i_k})  -  \sum_{k = 1}^{n-1} (-1)^{k+1} \sum_{2 \leq i_1 < \cdots < i_k \leq n} \gcd(d_1, d_{i_1}, \ldots, d_{i_k}) \\
&= \sum_{k = 1}^n (-1)^{k+1} \sum_{1 \leq i_1 < \cdots < i_k \leq n} \gcd(d_{i_1}, \ldots, d_{i_k}).
\end{align*}

We may now remove the restriction $\ord_p d_1 \geq \ord_p d_i$ as mentioned above. It is clear that
\[
\bb F_{q^{d_1}} \cap \bb F_{q^{d_2}} + \cdots + \bb F_{q^{d_1}} \cap \bb F_{q^{d_n}} \subseteq \bb F_{q^{d_1}} \cap (\bb F_{q^{d_2}} + \cdots + \bb F_{q^{d_n}}).
\]
By a similar dimension calculation made above, we see that these two $\bb F_q$-vector spaces have the same dimension, and thus they are equal.
\end{proof}

The formula for $e$ in (\ref{E: hyperplane}) quickly follows from the lemma since
\begin{align*}
N_1(\b d) &= \# \{ (x_1, \ldots, x_n) \in \bb F_{q^{d_1}} \times \cdots \times \bb F_{q^{d_n}} \mid a_1 x_1 + \cdots + a_n x_n = 0 \} \\
&= \dim_{\bb F_q} \bb F_{q^{d_1}} \cap (\bb F_{q^{d_2}} + \cdots + \bb F_{q^{d_n}}).
\end{align*}

In case it is of independent interest, we record here an alternative proof of Lemma \ref{L: Fq results} in the case $n = 3$.

\begin{proof}[Alternate proof]
Since $\bb F_{q^{d_1}} \cap \bb F_{q^{d_2}} + \bb F_{q^{d_1}} \cap \bb F_{q^{d_3}} \subseteq \bb F_{q^{d_1}} \cap (\bb F_{q^{d_2}} + \bb F_{q^{d_3}})$, if we show the dimensions of these spaces are equal, then the spaces must be equal. Let $m := \lcm(d_1, d_2, d_3)$. Let $\theta$ be a normal element of $\bb F_{q^m}$, and let $\sigma: x \mapsto x^q$ be the Frobenius morphism. Let $\alpha \in \bb F_{q^{d_2}}, \beta \in \bb F_{q^{d_2}}$, and $\gamma \in \bb F_{q^{d_1}}$ such that $\alpha + \beta = \gamma$. Using the normal basis,
\[
\alpha = \sum_{i=0}^{m-1} a_i \sigma^i(\theta) \in \bb F_{q^{d_2}} \Leftrightarrow \sigma^{d_2}(\alpha) = \alpha \Leftrightarrow a_i = a_{i + d_2} \text{ for all } i.
\]
Thus, we may identify $\alpha$ with the periodic sequence $(a_i)_{i=0}^\infty$ of period $d_2$, and note its generating function has the form
\[
G_\alpha := \sum_{i=0}^\infty a_i x^i = \frac{f_\alpha(x)}{1 - x^{d_2}} \qquad \text{with } \deg f_\alpha < d_2.
\]
Similarly, the generating functions for $\beta$ and $\gamma$ satisfy
\[
G_\beta(x) = \frac{f_\beta(x)}{1 - x^{d_3}}, \> \> G_\gamma(x) = \frac{f_\gamma(x)}{1 - x^{d_3}} \qquad \text{with } \deg f_\beta < m_3 \text{ and } \deg f_\gamma < d_1.
\]
Since $\alpha + \beta = \gamma$, we have
\[
\frac{f_\alpha(x)}{1 - x^{d_2}} + \frac{f_\beta(x)}{1 - x^{d_3}} = \frac{f_\gamma(x)}{1 - x^{d_1}},
\]
or
\begin{equation}\label{E: bunch o poly}
(1 - x^{d_1})(1 - x^{d_3}) f_\alpha(x) + (1 - x^{d_1})(1 - x^{d_2}) f_\beta(x) = (1 - x^{d_2})(1 - x^{d_3}) f_\gamma(x).
\end{equation}
Set $r := \gcd(d_2, d_3)$, and note that (\ref{E: bunch o poly}) shows
\[
(1 - x^{d_1}) \mid (1 - x^{d_2})(1 - x^{d_3})(1 - x^r)^{-1} f_\gamma(x).
\]
Thus, with
\[
h(x) := \gcd\left( (1 - x^{d_1}), (1 - x^{d_2})(1 - x^{d_3})(1 - x^r)^{-1} \right)
\]
we see that $(1 - x^{d_1})/h(x)$ divides $f_\gamma(x)$. That is, $f_\gamma(x) = g(x)(1 - x^{d_1})/h(x)$ for some $g$, and since $\deg f_\gamma < d_1$, we have $\deg g < \deg h$. Observe that this process is reversible: given any polynomial $g$ with $\deg g < \deg h$, if we set $f_\gamma(x) = g(x)(1 - x^{d_1})/h(x)$ then there exist polynomials $f_\alpha$ and $f_\beta$  of degrees $< d_2$ and $< d_3$, resp., satisfying (\ref{E: bunch o poly}). Since the set of such $g$ forms a vector space of dimension $\deg h$ over $\bb F_q$, we see that $\dim_{\bb F_q} \bb F_{q^{d_1}} \cap (\bb F_{q^{d_2}} + \bb F_{q^{d_3}}) = \deg h$. The result follows since 
\[
\deg \gcd\left( (1 - x^{d_1}), (1 - x^{d_2})(1 - x^{d_3})(1 - x^r)^{-1} \right) = \gcd(d_1, d_2) + \gcd(d_1, d_3) - \gcd(d_1, d_2, d_3).
\]
\end{proof}

\section{Example: Affine curve $y = x^n$}\label{S: affine curve}

In this section we will compute the partial zeta function of the affine curve $\c C_n := \{ y = x^n \}$ in $\bb A^2$. Set $c := \gcd(d_1, d_2)$, and denote by $\Phi_k$ the $k$-th cyclotomic polynomial over $\bb Q$. We will prove the following.

\begin{theorem}\label{T: curve}
\begin{equation}\label{E: curve product}
Z(\c C_n / \bb F_q, \b d, T) = \frac{1}{1-T} \cdot \prod_{k \mid \varphi(n)}\left( \frac{\Phi_k(T)}{\Phi_k(q^{c}T)} \right)^{a_k},
\end{equation}
where $a_k$ is an exponential sum defined by
\[
a_k := \frac{1}{\varphi(n)} \sum_{i = 1}^{\varphi(n)}  \gcd(n, M_i) \zeta_k^{i} \in \bb Z.
\]
Here $\varphi$ is Euler's totient function, $\zeta_k$ is a primitive $k$-th root of unity, and $M_i :=  (q_1^i - 1)/(q^{ci}-1)$.

Also, note that the total degree of $Z(\c C_n / \bb F_q, \b d, T)$ is bounded above by $1 + 2n^2$, independent of $q$ and $\b d$. This provides positive evidence for \cite[Question 2.5]{MR2030375}.
\end{theorem}

Before proving the theorem, we give a few remarks. First, exponential sums involving the greatest common divisor are well known. For example, $\varphi(n) = \sum_{r=1}^n \gcd(n, r) \zeta_n^{-r}$. Inspired by the above, we wonder:

\begin{conjecture}
Let $n, a$, and $d$ be positive integers. Set $M_i := (a^{di} - 1)/(a^i-1)$. We conjecture that the exponential sum $\frac{1}{\varphi(n)} \sum_{i = 1}^{\varphi(n)}  \gcd(n, M_i) \zeta_k^{i}$ is a rational integer. (Is there a closed form expression?)
\end{conjecture}

When $n$ is a prime number we give the following closed form for $a_k$. To state the result, define the Iverson bracket\footnote{If you are unfamiliar with the Iverson bracket, take a look at Donald Knuth's article \cite{MR1163629}} $[ \> ]$ for a logical proposition $\c P$ by:
\[
[\c P] := 
\begin{cases}
1 & \text{if $\c P$ is true} \\
0 & \text{if $\c P$ is false}.
\end{cases}
\]

\begin{theorem}
Suppose $n$ is a prime number. If $n = p$, then $a_k = [ k \mid 1]$. If $n \not= p$, then 
\[
a_k = \frac{\gcd(n, d_1/c) - n}{t} [k \mid t] + \frac{\phi(n)}{t_1} [k \mid t_1] + [k \mid 1],
\]
where $t := \ord_n^*(q^c) = $ the multiplicative order of $q^c$ modulo $n$, and $t_1 := \ord_n^*(q^{d_1})$.
\end{theorem}

\begin{proof}
If $n = p$ then $\gcd(n, M_i) = 1$, and so $a_k = [k \mid 1]$. Suppose now that $n$ is a prime number different from $p$. Observe that since $c \mid d_1$, we have $t_1 = t / \gcd(t, d_1/c)$. Now, as $M_i = 1 + q^{ci} + q^{2ci} + \cdots + q^{(d_1 - c)i}$, we have
\[
\gcd(n, M_i) = 
\begin{cases}
\gcd(n, d_1/c) & \text{if } t \mid i \\
\gcd(n, q^{d_1 i } - 1) & \text{otherwise}.
\end{cases}
\]
Moreover,
\[
\gcd(n, M_i) = 
\begin{cases}
\gcd(n, d_1/c) & \text{if } t \mid i \\
n & \text{if } t_1 \mid i \text{ and } t \nmid i \\
1 & \text{if } t_1 \nmid i.
\end{cases}
\]
Thus,
\begin{equation}\label{E: ak equ}
a_k = \frac{1}{\phi(n)}\left( \sum_{t \mid i} \gcd(n, d_1/c) \zeta_k^i + \sum_{t_1 \mid i, t \nmid i} n \zeta_k^i + \sum_{t_1 \nmid i} \zeta_k^i \right).
\end{equation}
For any divisor $t$ of $\phi(n)$ we have
\begin{equation}\label{E: divisor exp sum}
\sum_{\substack{i=1 \\ t \mid i}}^{\phi(n)} \zeta_k^i = \sum_{j=1}^{\phi(n)/t} \zeta_k^{tj} = \frac{\phi(n)}{t} [ k \mid t].
\end{equation}
Substituting (\ref{E: divisor exp sum}) into (\ref{E: ak equ}), then
\begin{align*}
a_k &= \frac{1}{\phi(n)} \left( \gcd(n, d_1/c) \frac{\phi(n)}{t} [k \mid t] + p \left( \frac{\phi(n)}{t_1} [ k \mid t_1] - \frac{\phi(n)}{t} [k \mid t] \right) + \phi(n) [k \mid 1] - \frac{\phi(n)}{t_1} k \mid t_1] \right) \\
&= \frac{\gcd(n, d_1/c) - n}{t} ;k \mid t] + \frac{\phi(n)}{t_1} [ k \mid t_1] + [k \mid 1]
\end{align*}
as desired.
\end{proof}

Next, let us illustrate Theorem \ref{T: curve} by explicitly computing the cases $n = 2$ and $3$. In the quadratic case $n = 2$,
\[
Z(\c C_2 / \bb F_q, \b d, T) = 
\begin{cases}
        \frac{1}{1- q^{c}T} & \text{if $q$ even; or $q$ odd and $d_1/c$ odd} \\
        \frac{ {1-T}}{{(1- q^{c}T)}^{2}} & \text{if $q$ odd and $d_1/c$ even.}
 \end{cases}
\]
For the cubic case $n = 3$, we split into residue classes. For $q \equiv$ $0$ or $1$ modulo $3$:
\[
Z(\c C_3 / \bb F_q, \b d, T) = 
\begin{cases}
\frac{1}{1-q^{c}T} & \text{if $q \equiv 0$ mod $3$} \\ 
& \text{if $q \equiv 1$ mod $3$ and $\frac{d_1}{c} \not\equiv 0$ mod $3$} \\
& \\
\frac{(1-T)^{2}}{( 1-q^cT )^{3}} & \text{if $q \equiv 1$ mod $3$ and $\frac{d_1}{c} \equiv 0$ mod $3$}
\end{cases}
\]
and for $q \equiv 2$ mod $3$:
\[
Z(\c C_3 / \bb F_q, \b d, T) = 
        \begin{cases}
            \frac{1}{1-q^{c}T} 
            & \text{\tiny $q \equiv 2\ \text{mod}\ 3, c\ \text{even}, \frac{d_1}{c} \not\equiv 0\ \text{mod}\ 3;$} \\
            & \text{\tiny $q \equiv 2\ \text{mod}\ 3, c\ \text{odd},\frac{d_1}{c}\ \text{odd}, \frac{d_1}{c} \not\equiv 0\ \text{mod}\ 3;$}\\ \\
            \frac{(1-T)^{2}}{\left( 1-q^{c}T\right)^{3}}
            & \text{\tiny $q \equiv 2\ \text{mod}\ 3, c\ \text{even}, \frac{d_1}{c} \equiv 0\ \text{mod}\ 3;$} \\
            & \text{\tiny $q \equiv 2\ \text{mod}\ 3, c\ \text{odd}, \frac{d_1}{c}\ \text{even},\frac{d_1}{c} \equiv 0\ \text{mod}\ 3; $}\\
            \\
            \frac{(1+q^{c}T)(1 - T)}{(1 - q^{c}T)^{2}(1+T)}
            & \text{\tiny $q \equiv 2\ \text{mod}\ 3, \frac{d_1}{c}\ \text{even}, \frac{d_1}{c} \not\equiv 0\ \text{mod}\ 3$}
            \\ \\
            \frac{(1-T)(1+T)}{(1-q^{c}T)^2(1+q^{c}T)}
            & \text{\tiny $q \equiv 2\ \text{mod}\ 3, c \ \text{odd}, \frac{d_1}{c} \equiv 0\ \text{mod}\ 3$}
        \end{cases}
 \]

We now move to the proof of Theorem \ref{T: curve}, which will consist of the rest of this section.

\begin{lemma}\label{L: count}
Let $\b d = (d_1, d_2) \in \bb Z_{\geq 1}^2$. Set $c := \gcd(d_1, d_2)$. For $\c C_n$ the affine curve $y = x^n$ in $\bb A^2$,
\[
N_m(\b d) = \gcd(n, M_m)(q^{c m} - 1) + 1
\]
where $M_m := (q_1^m - 1)/(q^{cm}-1) = \sum_{i = 0}^{(d_1/c) - 1}\, q^{cmi}$.
\end{lemma}

\begin{proof}
The number of solutions to $y = x^n$ is the number of elements in $\mathbb{F}_{q_1^m}$ whose $n$-th power lies in $\mathbb{F}_{q_2^m}$.  Let $\alpha$ be a generator of $\mathbb{F}_{q_1^m}^{\times}$. For notational convenience, for positive integers $i, j$ with $i$ dividing $j$, set
\[
[ i ,j ] := \{i, 2i, 3i, \ldots, j \}.
\] 
Then the number of nonzero solutions to $y = x^n$ is
\[
\# \{k \in [1, q_1^m - 1]\mid \alpha^{ nk } \in \mathbb{F}_{q_1^m} \cap \mathbb{F}_{q_2^m} \} = \# \{k \in  [1, q_1-1] \mid M_{m}\  \text{divides}\ nk\}.
\]
which follows since $M_{m}(q^{cm} - 1) = q_1^m - 1$ and $\alpha^{nk} \in \mathbb{F}_{q_1^m} \cap \mathbb{F}_{q_2^m}$ if and only if $\frac{q_1^m -1}{q^{cm} -1}$ divides $nk$. Therefore,
\begin{align*}
N_m(\b d) &= \# \{k \in [1 , q_1^m - 1] \mid M_{m}\ \text{divides}\ nk\} + 1  \\
&= \# \{k \in [n,n(q_1^m - 1)] \mid M_{m}\ \text{divides}\ k\} + 1   \\
&= \#\{k \in [1,n(q_1^m - 1)] \mid M_{m} \text{ and } n\ \text{divide}\ k\}+ 1 \\
&= \#\{k \in [M_{m},n(q^{cm} - 1)M_m] \mid n\ \text{divides}\ k\}+ 1.
\end{align*}
Now any $M_{m}k \in [M_m, n(q^{cm} - 1)M_{m}]$ is divisible by $n$ if and only if $n/\gcd(n, M_{m})$ divides $k$. Thus
\begin{align*}
\{M_{m}k : k \in [1, n(q^{cm} - 1)] \text{ and } n \mid Mk  \} &= \bigg\{ M_{m}k : k \in [1, n(q^{cm} -1)] \text{ and } \frac{n}{\gcd(n, M_{m})}\ \text{divides}\ k \bigg\} \\
                                                 &= \bigg\{ \frac{n M_{m}}{\gcd(n, M_{m})}k : k \in [1 ,\gcd(n, M_{m})(q^{cm} -1)]\bigg\},
\end{align*}
which means that
\begin{align*}
N_m(\b d) &=   \#\bigg\{ \frac{n M_{m}}{\gcd(n, M_{m})}k : k \in [1 ,\gcd(n, M_{m})(q^{cm} -1)]\bigg\}  + 1 \\
&=  \gcd(n, M)(q^{cm} -1) + 1.
\end{align*}
\end{proof}

\begin{lemma}\label{L: seq}
The sequence $\{ \gcd(n, M_m) \}_{m=1}^\infty$ has period $\varphi(n)$.
\end{lemma}

\begin{proof}
By Euler's theorem, $M_{m} \equiv M_{m'}$ mod $n$ whenever $m \equiv m'$ mod $\varphi(n)$. We arrive at the statement by noting that $\gcd(n, b) = \gcd(n, b \text{ mod } n)$, so in particular,
\[
\gcd(n, M_{m}) = \gcd(n, M_{m} \text{ mod } n) = \gcd(n, M_{m'} \text{ mod } n) = \gcd(n, M_{m'})
\]
whenever $m \equiv m'$ mod $\varphi(n)$.
\end{proof}

Applying Lemmas \ref{L: count} and \ref{L: seq},
\begin{align*}
Z(\c C_n / \bb F_q, \b d, T) &= \exp \left( \sum_{ i = 1}^{\varphi(n)}\sum_{\substack{m \equiv i \text{ mod } \varphi(n)  \\ m \geq 1}}\gcd(n, M_{i})(q^{mc} - 1)\frac{T^{m}}{m} + \frac{T^{m}}{m}\right) \\
     &= \exp \left( \sum_{ i = 1}^{\varphi(n)}\gcd(n, M_{i})\sum_{\substack{m \equiv i \text{ mod } \varphi(n)  \\ m \geq 1}}\frac{{(q^{c}T)}^{m}}{m} - \frac{T^{m}}{m}\right)/(1-T).
\end{align*}
Now using the fact that
\[
    \sum_{j = 1}^{\varphi(n)}\zeta^{(m-i)j} = 
    \begin{cases}
        \varphi(n) & \text{if $m \equiv i$ mod $\varphi(n)$} \\
        0 & \text{otherwise},
    \end{cases}
\]
where $\zeta$ is a primitive $\varphi(n)$-th root of unity, observe that 
\begin{align*}
 \sum_{\substack{m \equiv i \text{ mod } \varphi(n) \\ m \geq 1}} \frac{T^{m}}{m}
    &= -\frac{1}{\varphi(n)}\sum_{j = 1}^{\varphi(n)}\zeta^{-ij}\log(1-\zeta^{j}T).
\end{align*}
Set $\delta_{i} := \gcd(n, M_{i})$ and define the exponential sum
\[
    S_{j} := \sum_{i = 1}^{\varphi(n)} \delta_{i} \zeta^{-ij}.
\]
Then we may write
\begin{align}\label{E: zeta almost}
Z(\c C_n / \bb F_q, \b d, T) &= \exp \left( \sum_{ i = 1}^{\varphi(n)}\delta_{i}\sum_{\substack{m \equiv i \text{ mod } \varphi(n)  \\ m \geq 1}}\frac{{(q^{c}T)}^{m}}{m} - \frac{T^{m}}{m}\right)/(1-T) \notag \\
 &= \exp \left( \frac{1}{\varphi(n)} \sum_{ i = 1}^{\varphi(n)} \delta_{i} \sum_{j = 1}^{\varphi(n)}\zeta^{-ij} \log\left(\frac{1-\zeta^{j}T}{1-\zeta^{j}q^{c}T} \right) \right)/(1-T) \notag \\
  &= \exp \left( \log\prod_{j = 1}^{\varphi(n)}\left(\frac{1-\zeta^{j}T}{1-\zeta^{j}q^{c}T} \right)^{S_{j}/\varphi(n)} \right)/(1-T) \notag\\
    &=  \frac{1}{1-T} \cdot \prod_{j = 1}^{\varphi(n)}\left(\frac{1-\zeta^{j}T}{1-\zeta^{j}q^{c}T} \right)^{S_{j}/\varphi(n)}.
\end{align}

At this point, it is not clear that $Z(\c C_n / \bb F_q, \b d, T)$ is a rational function over $\bb Q$ due to the exponents $S_j$ being exponential sums, and unfortunately, we are unable to prove they are integers. However, Wan's rationality theorem \cite{MR2027782} tells us it is a rational function over $\bb Q$, and thus since the reciprocal zeros and poles are distinct it must be the case that $S_j$ is an integer divisible by $\varphi(n)$. 

Denote by $o(j)$ the additive order of $j$ modulo $\varphi(n)$. Then
\[
Z(\c C_n / \bb F_q, \b d, T) =  \frac{1}{1-T} \cdot \prod_{k \mid \varphi(n)} \prod_{\substack{1 \leq j \leq \varphi(n) \\ o(j) = k}}  \left(\frac{1-\zeta^{j}T}{1-\zeta^{j}q^{c}T} \right)^{S_{j}/\varphi(n)}.
\]
Let $\sigma$ by an automorphism of $\bb Q(\zeta)(T)$ over $\bb Q(T)$, and denote by $j_\sigma$ the integer such that $1 \leq j_\sigma \leq \varphi(n)$ and $\sigma(\zeta^j) = \zeta^{j_\sigma}$. Then
\begin{align*}
Z(\c C_n / \bb F_q, \b d, T) &= \sigma( Z(\c C_n / \bb F_q, \b d, T) ) \\
&= \frac{1}{1-T} \cdot \prod_{k \mid \varphi(n)} \prod_{\substack{1 \leq j \leq \varphi(n) \\ o(j) = k}}  \left(\frac{1- \zeta^{j_\sigma}T}{1-\zeta^{j_\sigma}q^{c}T} \right)^{S_{j}/\varphi(n)}.
\end{align*}
Since there always exists an automorphism $\sigma$ for which $j_\sigma = j'$ for any $j$ and $j'$ having additive orders $k$ mod $\varphi(n)$, we see that $S_{j'} = S_j$ for all $j, j'$ with additive order $k$ mod $\varphi(n)$, otherwise $\sigma( Z(\c C_n / \bb F_q, \b d, T) )$ would have zeros and poles with multiplicities different from $Z(\c C_n / \bb F_q, \b d, T)$. We may now define
\[
a_k := \frac{1}{\varphi(n)} \sum_{i = 1}^{\varphi(n)}  \gcd(n, M_i) \zeta^{-ij},
\]
where $j$ is any positive integer with additive order $k$ modulo $\varphi(n)$; note that $\zeta_k := \zeta^j$ is a primitive $k$-th root of unity, and so we may write
\begin{equation}\label{E: ak}
a_k := \frac{1}{\varphi(n)} \sum_{i = 1}^{\varphi(n)}  \gcd(n, M_i) \zeta_k^{i}.
\end{equation}
Finally,
\begin{align*}
Z(\c C_n / \bb F_q, \b d, T) &= \frac{1}{1-T} \cdot \prod_{k \mid \varphi(n)} \prod_{\substack{1 \leq j \leq \varphi(n) \\ o(j) = k}}  \left(\frac{1-\zeta^{j}T}{1-\zeta^{j}q^{c}T} \right)^{S_{j}/\varphi(n)} \\
&= \frac{1}{1-T} \cdot \prod_{k \mid \varphi(n)} \prod_{\substack{1 \leq j \leq \varphi(n) \\ o(j) = k}}  \left(\frac{1-\zeta^{j}T}{1-\zeta^{j}q^{c}T} \right)^{a_k} \\
&= \frac{1}{1-T} \cdot \prod_{k \mid \varphi(n)}  \left( \prod_{\substack{1 \leq j \leq \varphi(n) \\ o(j) = k}}  \frac{1-\zeta^{j}T}{1-\zeta^{j}q^{c}T} \right)^{a_k} \\
&= \frac{1}{1-T} \cdot \prod_{k \mid \varphi(n)}  \left(  \frac{\Phi_k(T)}{\Phi_k(q^c T)} \right)^{a_k},
\end{align*}
where $\Phi_k$ is the $k$-th cyclotomic polynomial. This prove \eqref{E: curve product}.

\section{Partial zeta function}\label{S: partial zeta}

In this section, we prove Theorem \ref{T: main zeta} and Corollary \ref{C: zeta ord}. The proof is a slight generalization of the argument given in \cite{MR1314464}; we provide the details for completeness. For each $i = 1, \ldots, n$, let $\{ \mu_{ij} \}_{j=1}^{a d_i}$ be a basis of $\bb F_{q^{d_i}} \subseteq \bb F_{q^d}$ over $\bb F_p$. Let $\{ \mu_j \}_{j=1}^{a d}$ be a basis of $\bb F_{q^d}$ over $\bb F_p$. We now make a change of variable: set $x_i := \sum_{j=1}^{a d_i} z_{ij} \mu_{ij}$. Observe that a monomial $x_1^{u_1} \cdots x_n^{u_n}$ transforms as follows. Writing each $u_i$ in base $p$, such as $u_1 = a_1 p^{\alpha_1} + \cdots$ and $u_2 = a_2 p^{\alpha_2} + \cdots $ with $0 \leq a_i \leq p-1$, we see that
\begin{align*}
x_1^{u_1} \cdots x_n^{u_n} &= \left( \sum_{j=1}^{a d_1} z_{1j} \mu_{1j} \right)^{u_1} \cdots \left(  \sum_{j=1}^{a d_n} z_{nj} \mu_{nj} \right)^{u_n} \\
&= \left( \sum_{j=1}^{a d_1} z_{1j} \mu_{1j}^{p^{\alpha_1}} \right)^{a_1} \cdots \left(  \sum_{j=1}^{a d_n} z_{nj} \mu_{nj}^{p^{\alpha_n}} \right)^{a_n} \cdots \\
&= \sum_{l=1}^{ad} G_l( z) \mu_l
\end{align*}
where each $G_l$ is a polynomial in the variables $z = (z_{ij})$, defined over $\bb F_p$. 

Performing this procedure for each $F_i$ creates a set of polynomials $\{ G_l^{(i)}(z) \}$ for $1 \leq l \leq ad$ and $1 \leq i \leq r$ such that $\deg \> G_l^{(i)} \leq w_p(F_i)$ and with the property
\[
N_1(\b d) = \# \{ z \in \bb F_p^{a(d_1 + \cdots + d_n)} \mid  G_l^{(i)}(z) = 0 \text{ for every $l$ and $i$} \}.
\]
That $p^\omega \mid N_1(\b d)$ now follows from the well-known theorem of Katz \cite{Katz-On_thm_of_Ax}. This proves Theorem \ref{T: main zeta}.

With
\[
\rho := \frac{(d_1 + \cdots + d_n) - d \sum_{i=1}^r w_p(F_i)}{\max_i w_p(F_i)},
\]
observe that $\ord_q N_m(\b d) \geq m \rho$ for every $m \geq 1$. Consequently, the argument of \cite[p. 256]{Ax-Zeros_of_polynomials_over_finite_fields} proves Corollary \ref{C: zeta ord}.

\section{Partial mixed character sums}\label{S: mixed}

Fix positive integers $d_1, \ldots, d_n \in \bb Z_{\geq 1}$, and set $d := \lcm(d_1, \ldots, d_n)$. Let $f(x) = \sum_{u \in U} a_u x^u \in \bb F_{q^d}[x_1, \ldots, x_n]$.  Let $\psi$ be a nontrivial additive character on $\bb F_{q^d}$, and $\chi_1, \ldots, \chi_r$ multiplicative characters on $\bb F_{q^{d_1}}^\times, \ldots, \bb F_{q^{d_n}}^\times$, respectively, trivial or non-trivial. Define the mixed partial character sum
\begin{align*}
S({\bf \chi}, {\bf d}, f) := \sum \chi_1(x_1) \cdots \chi_r(x_r) \psi( f(x) ),
\end{align*}
where the sum runs over $x_i \in \bb F_{q^{d_i}}^\times$ for $i = 1, \ldots, r$ and $x_i \in \bb F_{q^{d_i}}$  for  $i = r+1, \ldots, n$. Our main theorem of this section provides an estimate for $\ord_q S({\bf \chi}, {\bf d}, f)$. We first embed $S({\bf \chi}, {\bf d}, f)$ into the $p$-adic numbers. 

Let $\bb Q_p$ be the $p$-adic numbers, and denote by $\bb Q_{q^d}$ the unramified extension of $\bb Q_p$ of degree $ad$ (recall, $q = p^a$). We normalize the valuation such that $\ord_q(q) = 1$. Denote by $T_i$ the set of solutions of $z^{q_i} = z$ in $\bb Q_{q^{d_i}} \subseteq \bb Q_{q^d}$. Set
\[
\c T := T_1^\times \times \cdots \times T_r^\times \times T_{r+1} \times \cdots \times T_n.
\]
There exist $e_i \in \{1, \ldots, q_i - 1\}$ for $i = 1, \ldots, r$ such that  
\[
S(\chi, \b d, f) = \sum_{x \in \c T} x_1^{e_1} \cdots x_r^{e_r} \psi(f(\bar x)) \in \bb Q_{q^d}(\zeta_p),
\]
where $\bar x_i$ denotes the image of $x_i$ in the residue field $\bb F_{q^{d_i}} \subseteq \bb F_{q^d}$.

\begin{theorem}\label{T: ord S}
\[
\ord_q S(\chi, \b d, f) \geq \frac{1}{w_p(f)} \left( \sum_{i = 1}^n d_i -  \frac{1}{a(p-1)} \sum_{i=1}^r \sigma_p(e_i)  \right).
\]
When there is no twist (ie. $r = 0$), the estimate becomes
\[
\ord_q S(\b d, f) \geq (d_1 + \cdots + d_n) / w_p(f).
\]
\end{theorem}

\begin{proof}
Our proof follows that of \cite{MR932797} and \cite{MR2044049}, whose roots go back to at least \cite{Ax-Zeros_of_polynomials_over_finite_fields}. Denote by $P$ the polynomial $P(t) = \sum_{m=0}^{q^d -1} c_m t^m$ with the property that $P(z) = \psi(\bar z)$ for every $(q^d-1)$-root of unity $z$ in $\bb Q_{q^d}$. Note that the coefficients satisfy
\[
c_0 = 1, \qquad c_{q^d-1} = - q^d / (q^d-1), \qquad \text{and } c_m = g_m/(q^d-1) \text{ for } 1 \leq m \leq q^d-2,
\]
where $g_m$ is the Gauss sum
\[
g_m := \sum_{z^{q^d-1} = 1, z \in \bb Q_{q^d}} z^{-m} \psi(\bar z).
\]
By a well-known result of Stickelberger,
\begin{equation}\label{E: Stickelberger}
\text{for } 0\leq m \leq q^d-1: \qquad \ord_q c_m = \frac{\sigma_p(m)}{a(p-1)}.
\end{equation}
Recall, $f(x) = \sum_{u \in U} a_u x^u \in \bb F_{q^d}[x_1, \ldots, x_n]$. Let $\hat a_u \in \bb Q_{q^d}$ be the Teichm\"uller lift of $a_u$. Set $\b e := (e_1, \ldots, e_r, 0, \ldots, 0) \in \bb Z_{\geq 0}^n$. Then
\begin{align*}
S(\chi, \b d, f) &=  \sum_{x \in \c T} x^{\b e} \psi( f(\bar x) )  \\ 
&= \sum_{x \in \c T} x^{\b e} \prod_{u \in U} \psi( a_u \bar x^u )  \\
&= \sum_{x \in \c T} x^{\b e} \prod_{u \in U} P( \hat a_u  x^u ) \\
&= \sum_{x \in \c T} x^{\b e} \prod_{u \in U} \left( \sum_{m=0}^{q^d-1} c_m \hat a_u^m x^{m u} \right).
\end{align*}
Denote by $\Phi$ the set of functions $\phi: U \rightarrow \{ 0, 1, \ldots, q^d - 1 \}$. Then
\begin{align}\label{E: S expand}
S(\chi, \b d, f) &= \sum_{x \in \c T} x^{\b e} \sum_{\phi \in \Phi} \left( \prod_{u \in U} c_{\phi(u)} \hat a_u^{\phi(u)} \right) x^{\sum_{u \in U} \phi(u) u} \notag \\
&= \sum_{\phi \in \Phi} \left( \prod_{u \in U} c_{\phi(u)} \hat a_u^{\phi(u)} \right) \left( \sum_{x \in \c T} x^{\b e + \sum_{u \in U} \phi(u) u} \right).
\end{align}
The right most sum may be explicitly computed. In order to give its formula, we define the following. Denote by $m(\phi)$ the set of $i \in \{ r+1, \ldots, n\}$ such that the $i$-th entry of $\sum_{u \in U} \phi(u) u$ is non-zero, and $(m(\phi))$ the complement of $m(\phi)$ in $\{ r+1, \ldots, n\}$. Set
\[
q_{(m(\phi))} := \prod_{i \not\in m(\phi)} q_i \qquad \text{and} \qquad (q-1)_{m(\phi)} := \prod_{i \in m(\phi)} (q_i - 1).
\]
Set $\b{(q_* - 1)} := (q_1-1, \ldots, q_n-1)$. Also, define
\[
\b{(q_* - 1)}( \bb Z_{>0}^r \times \bb Z_{\geq 0}^{n-r}) := (q_1 - 1) \bb Z_{>0} \times \cdots \times (q_r - 1) \bb Z_{>0} \times (q_{r+1} - 1) \bb Z_{\geq 0} \times \cdots \times (q_n - 1) \bb Z_{\geq 0}.
\]
Then, using the well-known identities
\[
j \in \bb Z: \qquad \sum_{t \in T_i^\times} t^j = 
\begin{cases}
q_i - 1 & \text{if } (q_i - 1) \mid j \\
0 & \text{otherwise},
\end{cases}
\]
and
\[
j \in \bb Z_{\geq 0}: \qquad \sum_{t \in T_i} t^j = 
\begin{cases}
q_i - 1 & \text{if } (q_i - 1) \mid j \text{ and } j \not= 0 \\
q_i & \text{if } j = 0 \\
0 & \text{if } (q_i - 1) \nmid j,
\end{cases}
\]
we have
\begin{align*}
\sum_{x \in \c T} &x^{\b e + \sum_{u \in U} \phi(u) u}  \\
&=
\begin{cases}
q_{(m(\phi))} (q-1)_{m(\phi)} (q_1 - 1) \cdots (q_r - 1) & \text{if } \b e + \sum_{u \in U} \phi(u) u \in \b{(q_* - 1)}( \bb Z_{>0}^r \times \bb Z_{\geq 0}^{n-r}) \\
0 & \text{otherwise.}
\end{cases}
\end{align*}
Denote by $\Phi_0$ the set of functions $\phi: U \rightarrow \{0, 1, \ldots, q^d - 1\}$ such that $\b e + \sum_{u \in U} \phi(u) u \in \b{(q_* - 1)}( \bb Z_{>0}^r \times \bb Z_{\geq 0}^{n-r})$. 
It follows now by (\ref{E: S expand}) that
\[
\ord_q S(\chi, \b d, f) \geq \min_{\phi \in \Phi_0} \left\{ \sum_{i \not\in m(\phi)} d_i + \sum_{u \in U} \ord_q c_{\phi(u)} \right\}.
\]
Let $A \subseteq \{r+1, \ldots, n\}$, and denote by $\Phi_0^{(A)}$ the set of functions $\phi: U \rightarrow \{0, 1, \ldots, q^d - 1\}$ such that $\b e + \sum_{u \in U} \phi(u) u \in \b{(q_* - 1)}( \bb Z_{>0}^r \times \bb Z_{\geq 0}^{n-r})$ and the $i$-th entry of $\sum_{u \in U} \phi(u) u$ is non-zero if $i \in A$, and zero if $i \not\in A$. Then
\begin{align*}
\ord_q S(\chi, \b d, f) &\geq \min_{A \subseteq \{r+1, \ldots, n\}} \left\{ \sum_{i \not\in A} d_i + \min_{\phi \in \Phi_0^{(A)}} \left\{ \sum_{u \in U} \ord_q c_{\phi(u)} \right\} \right\} \\
&\geq \min_{A \subseteq \{r+1, \ldots, n\}} \left\{ \sum_{i \not\in A} d_i + \frac{1}{a(p-1)}\min_{\phi \in \Phi_0^{(A)}} \left\{ \sum_{u \in U} \sigma_p(\phi(u)) \right\} \right\},
\end{align*}
where we used (\ref{E: Stickelberger}) for the second inequality. 

For a vector $u = (u_1, \ldots, u_n) \in \bb Z_{\geq 0}^n$, define $\sigma_p(u) := (\sigma_p(u_1), \ldots, \sigma_p(u_n))$. Set $|u| := u_1 + \cdots + u_n$. Also, for $u, v \in \bb Z_{\geq 0}^n$, we write $u \geq v$ if $u_i \geq v_i$ for every $i$. Denote by $\b d_A$ the vector in $\bb Z_{\geq 0}^n$ such that the $i$-th entry is 0 if $i \not\in A$, and $d_i$ if $i \in A$. Last, we recall some properties of $\sigma_p$; see \cite[Proposition 11]{MR2044049}. For $a, b \in \bb Z_{\geq 0}$,
\begin{enumerate}
\item $\sigma_p(a+b) \leq \sigma_p(a) + \sigma_p(b)$;
\item $\sigma_p(ab) \leq \sigma_p(a) \sigma_p(b)$;
\item Let $q = p^f$. If $c$ is a positive multiple of $(q-1)$, then $\sigma_p(c) \geq \sigma_p(q-1) = f(p-1)$;
\end{enumerate}
Using this, we see that for $\phi \in \Phi_0^{(A)}$,
\[
\sum_{u \in U} \sigma_p(\phi(u)) \sigma_p(u) + \sigma_p(\b e) \geq a(p-1) \b d_A.
\]
Thus,
\[
\sum_{u \in U} \sigma_p(\phi(u)) |\sigma_p(u)| + |\sigma_p(\b e)| \geq a(p-1) |\b d_A| = a(p-1) \sum_{i \in A} d_i.
\]
Since $w_p(f) := \max_{u \in U}  |\sigma_p(u)| $,
\[
\sum_{u \in U} \sigma_p(\phi(u)) \geq \frac{1}{w_p(f)} \left( a(p-1) \sum_{i \in A} d_i -  |\sigma_p(\b e)| \right),
\] 
and so
\begin{align*}
\ord_q S(\chi, \b d, f) &\geq \min_{A \subseteq \{r+1, \ldots, n\}} \left\{ \sum_{i \not\in A} d_i + \frac{1}{w_p(f)} \sum_{i \in A} d_i -  \frac{|\sigma_p(\b e)|}{a(p-1)w_p(f)}  \right\} \\
&\geq \frac{1}{w_p(f)} \left( \sum_{i = 1}^n d_i -  \frac{|\sigma_p(\b e)|}{a(p-1)}  \right).
\end{align*}
\end{proof}

\bibliographystyle{amsplain}
\bibliography{../References/References.bib}

\end{document}